\documentclass[11pt,a4paper]{article}

\usepackage{amsmath,amsfonts,amssymb,latexsym,graphics,epsfig,url}
\usepackage{color}
\usepackage{amsthm}
\usepackage[english]{babel}
\usepackage{graphicx}
\usepackage{soul}
\usepackage[utf8]{inputenc}

\newtheorem{theorem}{Theorem}[section]

\newtheorem{lemma}[theorem]{Lemma}
\newtheorem{corollary}[theorem]{Corollary}

\newtheorem{conjecture}[theorem]{Conjecture}

\DeclareMathOperator{\spec}{sp}
\def\v{\mbox{\boldmath $v$}}

\def\vecv{\mbox{\boldmath $v$}}
\def\vecw{\mbox{\boldmath $w$}}

\def\w{\mbox{\boldmath $w$}}

\def\vec0{\mbox{\boldmath $0$}}

\def\A{\mbox{\boldmath $A$}}
\def\B{\mbox{\boldmath $B$}}

\def\D{\mbox{\boldmath $D$}}

\def\L{\mbox{\boldmath $L$}}

\def\X{\mbox{\boldmath $X$}}

\def\1{\mbox{\boldmath $1$}}

\newcommand\restr[2]{\ensuremath{\left.#1\right|_{#2}}}

\begin{document}
	
	\title{Some bounds on the Laplacian eigenvalues of token graphs
		\thanks{This research has been supported by
AGAUR from the Catalan Government under project 2021SGR00434 and MICINN from the Spanish Government under project PID2020-115442RB-I00.
The research of M. A. Fiol was also supported by a grant from the  Universitat Polit\`ecnica de Catalunya with references AGRUPS-2022 and AGRUPS-2023.}
	}
	\author{C. Dalf\'o$^a$,  M. A. Fiol$^b$, and A. Messegu\'e$^{c}$\\
		\\
		{\small $^a$Dept. de Matem\`atica, Universitat de Lleida, Igualada (Barcelona), Catalonia}\\
		{\small {\tt cristina.dalfo@udl.cat}}\\
		{\small $^{b}$Dept. de Matem\`atiques, Universitat Polit\`ecnica de Catalunya, Barcelona, Catalonia} \\
		{\small Barcelona Graduate School of Mathematics} \\
		{\small  Institut de Matem\`atiques de la UPC-BarcelonaTech (IMTech)}\\
		{\small {\tt miguel.angel.fiol@upc.edu} }\\
		{\small $^c$Dept. de Matem\`atica, Universitat de Lleida, Igualada (Barcelona), Catalonia}\\
		{\small Dept. de Ciències de la Computació, Universitat Polit\`ecnica de Catalunya, Barcelona, Catalonia} \\
		{\small {\tt visitant.arnau.messegue@udl.cat, arnau.messegue@upc.edu}}
	}
	
	\date{October 15, 2022}
	\maketitle	
	
\begin{abstract}
The $k$-token graph $F_k(G)$ of a graph $G$ on $n$ vertices is the graph whose vertices are the ${n\choose k}$ $k$-subsets of vertices from $G$, two of which being adjacent whenever their symmetric difference is a pair of adjacent vertices in $G$.

It is known that the algebraic connectivity (or
second Laplacian eigenvalue) of $F_k(G)$ equals the algebraic connectivity $\alpha(G)$ of $G$.

In this paper, we give some bounds on the (Laplacian) eigenvalues 
of a $k$-token graph (including the algebraic connectivity) 
in terms of the $h$-token graph, with $h\leq k$. 
For instance, we prove that if $\lambda$ is an eigenvalue of $F_k(G)$, but not of $G$, then
$$
\lambda\ge k\alpha(G)-k+1.
$$
As a consequence, we conclude that if $\alpha(G)\geq k$, then $\alpha(F_h(G))=\alpha(G)$ for every $h\le k$.
\end{abstract}

	\noindent{\em Keywords:} Token graph, Laplacian spectrum, Algebraic connectivity, Binomial matrix.
	
	\noindent{\em MSC2010:} 05C15, 05C10, 05C50.
	
	%%%%%%%%%%%%%%%%%%%%%%%%%%%%%%%%%%%%%%%%%%%%%%%%%%%%%%%%%%%%%%%%%%%%%%%%%%%%%%%%%%%
	\section{Introduction}
	\label{sec:-1}
\begin{figure}[t]
		\begin{center}
			\includegraphics[width=15cm]{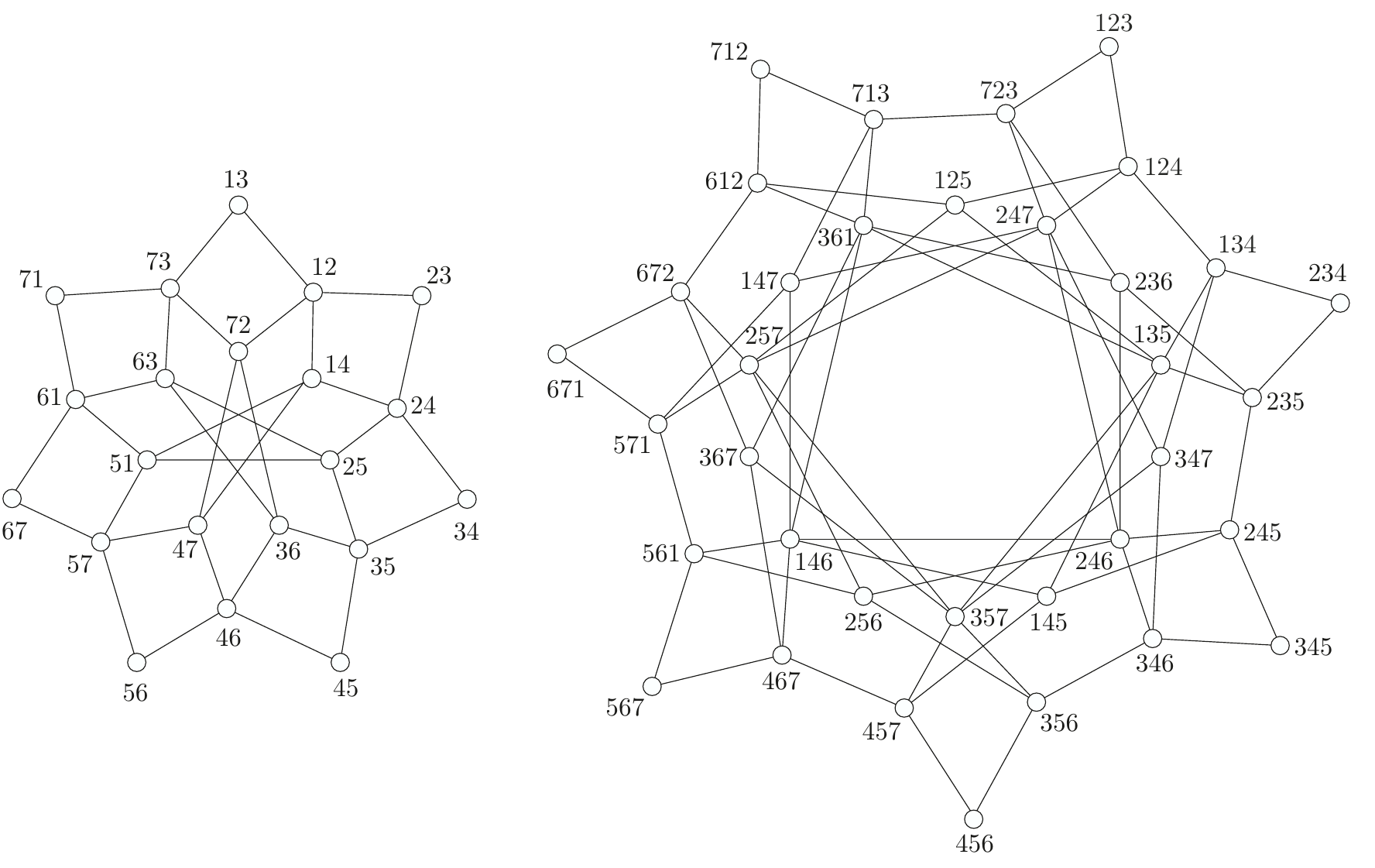}
			%\vskip -1.5cm
			\caption{The $2$-token graph $F_2(C_7)$ of the cycle graph $C_7$ (left), and the $3$-token graph $F_3(C_7)$ (right).}
			\label{fig1}
		\end{center}
	\end{figure}

	Let $G$ be a simple graph with vertex set $V(G)=[n]=\{1,2,\ldots,n\}$ and edge set $E(G)$. Let $\Delta(G)$ denote the maximum degree of $G$. For a given integer $k$ such that
	$1\le k \le n$, the {\em $k$-token graph} $F_k(G)$ of $G$ is the graph whose vertex set $V (F_k(G))$ consists of the ${n \choose k}$
	$k$-subsets of vertices of $G$, and two vertices $A$ and $B$
	of $F_k(G)$ are adjacent whenever their symmetric difference $A \bigtriangleup B$ is a pair $(a,b)$ such that $a\in A$, $b\in B$, and $\{a,b\}\in E(G)$.
	The name `token graph'
	was given by
	Fabila-Monroy,  Flores-Pe\~{n}aloza,  Huemer,  Hurtado,  Urrutia, and  Wood \cite{ffhhuw12}, because the vertices of $F_k(G)$ correspond to configurations
	of $k$ indistinguishable tokens placed at distinct vertices of $G$, where
	two configurations are adjacent whenever one configuration can be reached
	from the other by moving one token along an edge from its current position
	to an unoccupied vertex. %Thus,
%	the maximum degree of $F_k(G)$ satisfies
%	\begin{equation}
%		\label{DeltaFk}
%		\Delta(F_k(G))\le k\Delta(G).
%	\end{equation}
	In Figure \ref{fig1}, we show the 2-token graph and 3-token graph of the cycle $C_7$ on 7 vertices.
 
	Note that if $k=1$, then $F_1(G)\cong G$; and if $G$ is the complete graph $K_n$, then $F_k(K_n)\cong J(n,k)$, where $J(n,k)$ denotes the Johnson graph.
	
Token graphs have some applications in physics. For instance, a
relationship between token graphs and the exchange of Hamiltonian operators in
quantum mechanics is given in Audenaert, Godsil, Royle, and Rudolph \cite{agrr07}.
	
Let $G$ have Laplacian matrix $\L=\L(G)$, with eigenvalues $\lambda_1(=0),\lambda_2,\ldots,\lambda_n$. The algebraic connectivity of $G$ is $\alpha(G)=\lambda_2$. 
% \textcolor{red}{Recently, it was conjectured by Dalf\'o, Duque, Fabila-Monroy, Fiol, Huemer, Trujillo-Negrete, and Zaragoza Mart\'{\i}nez \cite{ddffhtz21} that the algebraic connectivity of $F_k(G)$ equals $\alpha(G)$.} 
For more information on algebraic connectivity, see the survey by Abreu \cite{abreu07}.
	
This paper is structured as follows. In Section \ref{sec:-7}, we present some of the known results on the algebraic connectivity. In Section \ref{sec:new}, we give some bounds on the eigenvalues of a $k$-token graph $F_k(G)$ (including the algebraic connectivity), the first one in terms of the algebraic connectivity $\alpha(G)$ of graph $G$. As a corollary, we prove that if $\alpha(G)\geq k$, then $\alpha(F_h(G))=\alpha(G)$ for every $h\le k$. Thus, we show that  $\alpha(F_k(G))=\alpha(G)$ holds if $\alpha(G)\geq n/2$. Finally, in Subsection \ref{sec:alternative}, we prove a general result involving some eigenvalues of $F_k(G)$ and $F_{k-1}(G)$ and the maximum degree of $G$.
	
	\section{Known results}
	\label{sec:-7}
	
 In the proof of our main result, we use the following concepts and lemmas.
	
	Given a graph $G=(V,E)$ of order $n=|V|$,
	we say that a vector $\v\in \mathbb{R}^n$ is an \textit{embedding} of $G$ if $\v^{\top}\1=0$, where $\1$ is the all-$1$ vector.
	Note that if $\v$ is a $\lambda$-eigenvector of $G$, with Laplacian eigenvalue $\lambda>0$, then it is an embedding of $G$.

\begin{lemma}[Rayleigh quotient]
	\label{lemm:Rayleigh}
	For a graph $G=(V,E)$ with Laplacian matrix $\L=\L(G)$, and  $\v\in V$, let
\begin{equation}
\label{Rayleigh}
\lambda_G(\v):=\frac{\v^{\top}\L\v}{\v^{\top}\v}=\frac{\sum\limits_{(x,y)\in E}[\v(x)-\v(y)]^2}{\sum\limits_{x\in V}\v^2(x)},
\end{equation}	
where
$\v(x)$ denotes the entry of $\v$ corresponding to the vertex $x\in V(G)$.
If $\v$ is an eigenvector of $G$, then its corresponding eigenvalue is $\lambda_G(\v)$.
Moreover, for an embedding $\v$ of $G$, we have
	\begin{equation}
		\label{bound-lambda(v)}
		\lambda_G(\v)\ge \alpha(G),
	\end{equation}
and we have equality  when $\v$ is an $\alpha(G)$-eigenvector of $G$.
	\end{lemma}

The following results and subsequent conjecture were given by Dalf\'o, Duque, Fabila-Monroy, Fiol, Huemer, Trujillo-Negrete, and Zaragoza Mart\'{\i}nez in \cite{ddffhtz21}.

\begin{lemma}[\cite{ddffhtz21}]
\label{lem:1}
\begin{itemize}
\item[$(i)$]
Let $G$ be a graph on $n$ vertices.
The (Laplacian) spectra of the token graphs $F_k(G)$ for $k=1,\ldots,\lfloor n/2\rfloor$ satisfy
\begin{equation}
\label{eq:inclusion}
\spec(G)=\spec(F_1(G))\subset \spec(F_2(G))\subset \cdots\subset \spec(F_{\lfloor n/2\rfloor}(G))
\end{equation}
and, hence,
\begin{equation}
\label{ineq-alphas}
\alpha(G)\ge \alpha(F_2(G))\ge \cdots \ge \alpha(F_{\lfloor n/2\rfloor}(G)).
\end{equation}
\item[$(ii)$] 
Let $\B$ the $(n,k)$-binomial matrix with rows indexed by the $k$-subsets $X\subset [n]$ and columns indexed by the elements  $x\in [n]$, with entries $\B_{Xx}=1$ if $x\in X$, and $\B_{Xx}=0$ otherwise.
If $\vecv$ is a $\lambda$-eigenvector of $F_k(G)$ and $\B^{\top}\vecv \neq \vec0$, then $\B^{\top}\vecv$ is a $\lambda$-eigenvector of $G$. Hence,
if $\lambda=\alpha(F_k(G))<\alpha(G)$, then $\B^{\top}\v=\vec0$.
\item[$(iii)$]
%Let us show a result by Dalf\'o and Fiol \cite{df22}.
%\begin{lemma}[\cite{df22}]
%\label{lem:eigenvectors}
%Let $G$ be a graph with $k$-token graph $F_k(G)$.
Given a subset $U\in V(G)$ with $|U|\le k-1$, let $S_U:=\{X\in V(F_k(G)):U\subset X\}$.
Let $H_U\cong F_{k-1}(G-U)$  be the subgraph of $F_k(G)$ induced by $S_U$. 
Let $\v$ be an eigenvector of $F_k(G)$ such that $\B^{\top}\v=\vec0$, and
$
\w_U:=\restr{\v}{S_U}
$.
Then, $\w_U$ is an embedding of $H_U$.
%\end{lemma}
\end{itemize}
\end{lemma}

	%\begin{proof}
	%Assume that the matrix $\B^{\top}$ has the first row indexed by   $a\in V(G)$. Then, we have
	%$$
	%\vec0=
	%\B^{\top}\v=
	%\left(
	%\begin{array}{c|c}
	%\1^{\top} & \vec0^{\top}\\
	%\hline
	%\B_1 & \B_2
	%\end{array}
	%\right)
	%\left(
	%\begin{array}{c}
	%\w_a \\
	%\hline
	%\w'_a
	%\end{array}
	%\right)=
	%\left(
	%\begin{array}{c}
	%\1^{\top}\w_a \\
	%\hline
	%\B_1\w_a+\B_2\w'_a
	%\end{array}
	%\right),
	%$$
	%\textcolor{red}{
		%$$
		%\vec0=
		%\B^{\top}\v=
		%\left(
		%\begin{array}{cc}
		%	\1^{\top} & \vec0^{\top}\\
		%	\B_1 & \B_2
		%\end{array}
		%\right)
		%\left(
		%\begin{array}{c}
		%	\w_a \\
		%	\w'_a
		%\end{array}
		%\right)=
		%\left(
		%\begin{array}{c}
		%	\1^{\top}\w_a \\
		%	\B_1\w_a+\B_2\w'_a
		%\end{array}
		%\right),
		%$$}
	%where $\1^{\top}$ is a row ${n-1\choose k-1}$-vector, $\vec0$ is a row  ${n-1\choose k}$-vector, $\B_1=\B(n-1,k-1)^{\top}$, and  $\B_2=\B(n-1,k)^{\top}$.
	%Then, $\1^{\top}\w_a=0$, so that $\w_a$ is an embedding of $H_a$. Furthermore, since $\v$ is an embedding of $G$, we have $\1^{\top}\v=\1^{\top}\w_a+\1^{\top}\w'_a=0$ (with the appropriate dimensions of the all-1 vectors). Hence, it must be  $\1^{\top}\w'_a =0$, and $\w'_a$ is an embedding of $H'_a$.
	%\end{proof}
	
% 	The following conjecture was given by Dalf\'o, Duque, %Fabila-Monroy, Fiol, Huemer, Trujillo-Negrete, and Zaragoza %Mart\'{\i}nez \cite{ddffhtz21}.
	
		\begin{conjecture}[\cite{ddffhtz21}]
		\label{conjecture}
		Let $G$ be a graph on $n$ vertices. Then, for every $k=1,\ldots,n-1$, the algebraic connectivity of the token graph $F_k(G)$ equals the algebraic connectivity of $G$.
		%$$
		%\alpha(F_k(G))=\alpha(G)\quad \mbox{for every $k=1,\ldots,|V|-1$}.
		%$$
	\end{conjecture}

 In different papers  (see \cite{ddffhtz21,df22,rdfm23})
 it was proved that the conjecture holds for different infinite families of graphs. For instance,
 in \cite{df22}, the first two authors proved the following result.
\begin{theorem}[\cite{df22}]
	\label{theo:Delta}
	Let $G$ be a graph on $n$ vertices 
	%satisfying $\alpha(F_{k-1}(G))=\alpha(G)$ 
	and with minimum degree
	\begin{equation}
		\label{boundDelta}
		\delta(G)\ge %\phi(k)  =
  \frac{k(n+k-3)}{2k-1}
	\end{equation}
	for some integer $k\in\{1,\ldots, \lfloor n/2\rfloor\}$. Then,   the algebraic connectivity of its $h$-token graph equals the algebraic connectivity of $G$,
	$$
	\alpha(F_h(G))=\alpha(G),
	$$
	for every $h\le k$.
\end{theorem}
 
Actually, after submitting the first version of this paper, the authors learned (from Fabila-Monroy \cite{f23}) that this conjecture was already known as the {\em  Aldous' spectral gap conjecture},  and it was proved in 2010 by Caputo, Ligget, and Richthammer in \cite{clr10}.
 % Moreover,  Ouyang  \cite{o19} and Lew \cite{l23}
 % also mentioned that this conjecture was actually solved.
 Moreover, Cesi \cite{c16} provided a simpler proof of the so-called `octopus inequality', which is one of the main ingredients to prove the Aldous’
conjecture.  
These results were
obtained by using the theory of stochastic processes 
(continuous Markov chains of random walks and the 
interchange process). Thus, we think that our algebraic and combinatorial method, although not completely solving the conjecture, provides insight into the other eigenvalues of the token graphs.

% We use the following results by Fiedler \cite{fi73}.

% \begin{lemma}[\cite{fi73}]
% 	\textcolor{red}{If $G_1$ and $G_2$ are edge-disjoint graphs with the same set of vertices, then %$\alpha(G_1)+\alpha(G_2)=\alpha(G_1\cup G_2)$}.
% \end{lemma}

%  \begin{lemma}[\cite{fi73}]
% 	\textcolor{red}{Given a graph $G$ and its complementary graph $\overline{G}$ on $n$ vertices, then $\alpha(G)+\alpha(\overline{G})\leq n$, with equality if and only if $G$ is a complete graph or a null graph.}
% \end{lemma}

\section{New results}
\label{sec:new}

Here we present the main results of this paper, in which the more general of them is Theorem \ref{th:arnau}.
% The main results of this paper are the following theorem, where the latter is a slight improvement of the former.
\begin{theorem}
\label{TH:Fk-alpha}
If $\alpha(G)>\alpha(F_2(G))>\cdots >\alpha(F_k(G))$, 
we would have
%the algebraic connectivities of $F_k(G)$ and $G$ satisfy the %following inequality
\begin{equation}
\label{eq:Fk-alpha}
\alpha(F_k(G))\ge k[\alpha(G)-k+1].
\end{equation}
\end{theorem}
	
\begin{proof}
To apply induction, we first prove the result for $k=2$.
		Let $\v$ be an eigenvector of $F_2=F_2(G)$ with eigenvalue $\alpha(F_2(G))$ and norm $\|\v\|=1$.
		As before, let $S_x:=\{X\in V(F_2):x\in X \}$.
		Let $H_x\cong F_{1}(G-x)\cong G-x$ be the subgraph of $F_2$ induced by $S_x$.
		Let $\w_x:=\restr{\v}{H_x}$. Now, from the hypothesis $\alpha(G)>\alpha(F_2(G))$ and  Lemma \ref{lem:1}$(ii)$-$(iii)$, we have that $\w_x$ is an embedding of $H_x$. Hence,
		\[
		\lambda(\w_x)=\frac{\sum\limits_{(X,Y)\in E(H_x)} [\w_x(X)-\w_x(Y)]^2}{\sum\limits_{X\in V(H_x)}\w_x(X)^2}\ge \alpha(G-x)\ge \alpha(G)-1,
		\]
where the last inequality is due to Fiedler \cite{fi73}.
Then, using this,
		\begin{align}
			\alpha(F_2) =\lambda(\v)&=\sum\limits_{(X,Y)\in E(F_2)}[\v(X)-\v(Y)]^2
=\sum_{x\in V} \sum\limits_{(X,Y)\in E(H_x)}[\w_x(X)-\w_x(Y)]^2  \label{eq:vt-1}\\
			& \ge (\alpha(G)-1)\sum_{x\in V}\sum\limits_{X\in V(H_x)}\v(X)^2  =2[\alpha(G)-1].
			\label{eq:vt-2}
		\end{align}
		The last equality in \eqref{eq:vt-1} holds since, in the double summation, each edge $\{X,Y\}$ of $F_2$ is considered once, whereas the equality in \eqref{eq:vt-2} is because each vertex $X$ of $F_2$ is considered $k=2$ times.
Now, assume that the result holds for $k-1$ with $k\ge 3$, so that $\alpha(F_{k-1}(G))\ge (k-1)[\alpha(G)-k+2]$. Then, using the same notation as above with $H_x\cong F_{k-1}(G-x)$ (since $\alpha(F_k(G))<\alpha(F_{k-1}(G))$, the vector $\vecw_x$ is an embedding of $H_x$),
\begin{align*}
\lambda(\w_x)& =\frac{\sum\limits_{(X,Y)\in E(H_x)} [\w_x(X)-\w_x(Y)]^2}{\sum\limits_{X\in V(H_x)}\w_x(X)^2}
\ge \alpha(F_{k-1}(G-x))\\
 &\ge (k-1)[\alpha(G-x)-k+2]\ge (k-1)[\alpha(G)-k+1],
 \end{align*}
 where we used the induction hypothesis with $G-x$.
Hence, \eqref{eq:vt-1} and \eqref{eq:vt-2} become
\begin{align*}
\alpha(F_k) =\lambda(\v)&=\sum\limits_{(X,Y)\in E(F_k)}[\v(X)-\v(Y)]^2
=\frac{1}{k-1}\sum_{x\in V} \sum\limits_{(X,Y)\in E(H_x)}[\w_x(X)-\w_x(Y)]^2  \\
			& \ge \frac{1}{k-1}(k-1)[\alpha(G)-k+1]\sum_{x\in V}\sum\limits_{X\in V(H_x)}\v(X)^2=k[\alpha(G)-k+1],
			\end{align*}
since,  in the double summation, each edge $\{X,Y\}$ of $F_k$ is considered $k-1$ times, whereas for the last equality, each vertex $X$ of $F_k$ is considered $k$ times.
\end{proof}

\begin{corollary}
\label{coro:alpha>=k}
If $\alpha(G)\ge k$, then $\alpha(F_h(G))=\alpha(G)$ for every $h\le k$.
\end{corollary}
\begin{proof}
The proof is by contradiction. If $h\le k\le \alpha(G)$ and $\alpha(G)>\alpha(F_h(G))$, Theorem \ref{TH:Fk-alpha} yields $\alpha(G)> (h-1)[\alpha(G)-h+2]$, whence $\alpha(G)<h$ contradicting the hypothesis. Thus, it must be $\alpha(G)\le \alpha(F_h(G))$, which (together with \eqref{ineq-alphas}) allows us to conclude that $\alpha(F_h(G))=\alpha(G)$.
\end{proof}

% Theorem \ref{TH:Fk-alpha} also provides a relationship between the algebraic connectivities of the $k$-token and $h$-token of $G$.

% \begin{corollary}
% Let $h\le k$. Under the same hypotheses in Theorem \ref{TH:Fk-alpha}, we would have
% \begin{equation}
% \label{th:Fk-alpha}
% \alpha(F_k(G))\ge k\left[\frac{\alpha(F_h(G))}{h}-k+h\right].
% \end{equation}
% \end{corollary}	
% \begin{proof}
% Get an upper bound for  $\alpha(G)$ from $\alpha(F_h(G))\ge h[\alpha(G)-h+1]$ and substitute it into \eqref{eq:Fk-alpha}.
% \end{proof}

% In particular, if $h=k-1$, we would get 
% \begin{equation}
% \label{eq:alphaFk}
% \alpha(F_k(G))\ge k\left[\frac{\alpha(F_{k-1}(G)}{k-1}-1\right], 
% \end{equation}
% as we will prove in the following subsection by using another method.

Let us show some examples of graphs satisfying Corollary \ref{coro:alpha>=k}. The Hamming graph $H(d,q)$ has vertex set $[q]^d$, the set of ordered $d$-tuples of elements of $[q]$, or sequences of length $d$ from $[q]$,  and two vertices are adjacent if they differ in precisely one coordinate. In particular, $H(d,2)=Q_d$, the $d$-cube.
Then, since $\alpha(H(d,q))=q$, we have the following result.
\begin{corollary}
$\alpha(F_h(H(d,q)))=\alpha(H(d,q))=q$ for every $h\le q$.
\end{corollary}

In fact, since the Hamming graph is $H(d,q)=K_q\times\stackrel{(d)}{\cdots}\times K_q$, the Cartesian product of $d$ copies of the complete graph $K_q$, the last corollary is a direct consequence of the following result. 

\begin{lemma}
    Let $G_1$ and $G_2$ be two graphs such that $k\leq \alpha(G_1)\leq \alpha(G_2)$. Then, $\alpha(F_h(G_1\times G_2))=\alpha(G_1\times G_2)$ for every $h\leq k$.
\end{lemma}
\begin{proof}
    The result follows from Corollary \ref{coro:alpha>=k} and the result by Fiedler \cite{fi73} stating that
    $\alpha(G_1\times G_2)=\min\{\alpha(G_1),\alpha(G_2)\}$.
\end{proof}

Since we already know that the algebraic connectivity of $F_k(G)$ coincides with that of $G$, it is more interesting to use the method in the proof of Theorem \ref{TH:Fk-alpha} to give bounds on other eigenvalues of  $F_k(G)$. In fact, basically the same proof holds when, instead of $\alpha(F_k(G))$, we consider an eigenvalue $\lambda$ in $F_k(G)$ that is not in $G$. Thus, to ensure that we are dealing with embeddings of the involved subgraphs (implied by $\B^{\top}\v=0$) we can rid of the previous hypothesis $\alpha(G) > \alpha(F_2(G)) > \cdots > \alpha(F_k(G))$.
In the following theorem, we follow this approach, and just for completeness, we detail again the corresponding steps of the proof.

\begin{theorem}
\label{TH:Fk-alpha(b)}
Let $\lambda$ be an eigenvalue of $F_{k}=F_k(G)$ not in $G$. Then,
\begin{equation}
\label{eq:Fk-alpha(b)}
\lambda\ge k[\alpha(G)-k+1].
\end{equation}
\end{theorem}

\begin{proof}
		Let $\v$ be an eigenvector of $F_k(G)$ with eigenvalue $\lambda$ and norm $\|\v\|=1$.
		For any subset $U\in V(F_{k-1})$, let $S_U:=\{X\in V(F_k):U\in X \}$.
		Let $H_U\cong F_{1}(G-U)\cong G-U$ be the subgraph of $F_k$ induced by $S_U$.
		Let $\w_U:=\restr{\v}{H_U}$. Now, from the hypothesis and  Lemma \ref{lem:1}$(iii)$, we have that $\w_U$ is an embedding of $H_U$. Hence,
		\[
		\lambda(\w_U)=\frac{\sum\limits_{(X,Y)\in E(H_U)} [\w_U(X)-\w_U(Y)]^2}{\sum\limits_{U\in V(H_U)}\w_U(X)^2}\ge \alpha(G-U)\ge \alpha(G)-k+1,
		\]
where the last inequality is due to Fiedler \cite{fi73}.
Then, using this,
		\begin{align}
			\lambda=\lambda(\v)&=\sum\limits_{(X,Y)\in E(F_k)}[\v(X)-\v(Y)]^2
=\sum_{U\in V(F_{k-1})} \sum\limits_{(X,Y)\in E(H_U)}[\w_U(X)-\w_U(Y)]^2  \label{eq:vt-1(b)}\\
			& \ge (\alpha(G)-k+1)\sum_{U\in V(F_{k-1})}\sum\limits_{X\in V(H_U)}\v(X)^2  =2[\alpha(G)-k+1].
			\label{eq:vt-2(b)}
		\end{align}
		The last equality in \eqref{eq:vt-1(b)} holds since, in the double summation, each edge $\{X,Y\}$ of $F_k$ is considered once, whereas the equality in \eqref{eq:vt-2(b)} is because each vertex $X$ of $F_k$ is considered $k$ times.
  \end{proof}

% Reasoning as above, we have that the spectral radius of $F_k(G)$, denoted by $\lambda_{\max}(F_k(G))$, satisfies the following result.
% \begin{lemma}
% Let $G=(V,E)$ be  a vertex-transitive graph with a vertex $i\in [n]$. Then, the following inequality holds:
% $$
% \lambda_{\max}(F_k(G))\le \frac{k}{k-1}\lambda_{\max} (F_{k-1}(G-i)).
% $$
% \end{lemma}

	%\begin{proof}
	%Let $\v$ be an eigenvector of $F_k=F_k(G)$ with eigenvalue $\lambda_{\max}(F_k(G))$ and norm $\|\v\|=1$.  As before, let $S_i:=\{A\in V(F_k):i\in A \}$. Let $H_i\cong F_{k-1}(G-i)$ be the subgraph of $F_k$ induced by $S_i$. Note that, since $G$ is vertex-transitive, the graph $H_i$ is independent of the vertex $i$. Let $\w_i:=\restr{\v}{S_i}$. Now, we have
	%\[
	%\lambda(\w_i)=\frac{\sum\limits_{(A,B)\in E(H_i)} [\w_i(A)-\w_i(B)]^2}{\sum\limits_{A\in V(H_i)}\w_i(A)^2}\leq \lambda_{\max}(H_i).
	%\]
	%and, using this,
	%\begin{align}
	%\lambda_{\max}(F_k) =\lambda(\v)&=\sum\limits_{(A,B)\in E(F_k)}[\v(A)-\v(B)]^2 =\frac{1}{k-1}\sum_{i=1}^n
	% \sum\limits_{(A,B)\in E(H_i)}[\w_i(A)-\w_i(B)]^2  \label{eq:vt-1}\\
	%	 & \leq \frac{1}{k-1}\lambda_{\max}(H_i)\sum_{i=1}^n\sum\limits_{A\in V(H_i)}\v(A)^2  =\frac{k}{k-1}\lambda_{\max}(H_i).
	%\label{eq:vt-2}
	% \end{align}
%The second equality in \eqref{eq:vt-1}) holds since, in the double summation, each edge $(A,B)$ of $F_k$ is considered $k-1$ times, whereas the last equality is because each vertex $A$ of $F_k$ is considered $k$ times.
%\end{proof}

\subsection{A general result}
\label{sec:alternative}

In this subsection, we follow an alternative method that
leads to a more general result containing  Theorem \ref{TH:Fk-alpha(b)} with $\lambda=\alpha(F_k(G))$.
In fact, although in the next theorem we reason with the algebraic connectivities of $F_k=F_k(G)$ and $F_{k-1}=F_{k-1}(G)$, its proof also works if we replace $\alpha(F_k)$ by an eigenvalue $\lambda$, that is in $F_k$, but not in $F_{k-1}$.

% In fact, the method leads to an improvement of such a theorem, with a result that also depends on the maximum degree of $G$. 

Define $[n]_k$ to be the collection of subsets from $[n] = \left\{1,2,\ldots,n\right\}$ having size $k$,  $[n\setminus \left\{z\right\}]_k$ the collection of subsets from $[n] \setminus \left\{z\right\}$ having size $k$,  and $\B$  the $(n,k)$-binominal matrix.
Let $\vecv$ be any vector, indexed by the elements of $[n]_k$, associated with the token graph $F_k(G)$. 
Let $X\in [n]_{k}$. Denote by $X_{-x} = X \setminus \left\{x\right\}$ and $X_{-x}^{+y}=\X_{-x} \cup \left\{ y \right\}$.
Similarly, for any subset $Z \in [n\setminus \{z\}]_k$, let $Z^{+z}$ denote  
$Z\cup \{z\}$.

\begin{theorem}
\label{th:arnau}
%If $\lambda = \alpha(F_k(G))$ is not an eigenvalue of $\L_1$, then
%Let $F_k=F_k(G)$
If $\alpha(F_{k}(G))$ is not an eigenvalue of $F_{k-1}(G)$, then the following holds:
\begin{equation}
\alpha(F_k(G)) \geq 
\left\{
\begin{array}{ll}
      \displaystyle\frac{k}{k-1}\alpha(F_{k-1}(G))-\frac{k}{k-2}\min\{k-2,\Delta(G)\} & \mbox{if } k>2, \\[.2cm]
    2\alpha(G)-2 & \mbox{if } k=2,
\end{array}
\right.
\label{eq:arnau}
\end{equation}
where $\Delta(G)$ is the maximum degree of $G$.
\end{theorem}

\begin{proof}
  By the definition of the $k$-token graph of $G=(V,E)$, the degree $d_X$ of vertex $X$ in $F_k(G)$ is 
 $$
 d_X=\sum_{x\in X} d_x-\sum_{x,y \in X}a_{xy},
 $$
 where $d_x$ is the degree of $x$ in $G$, and $a_{xy}$ is the $xy$-entry of the adjacency matrix $\A$ of $G$.
 Moreover, by the same definition,  the entry of the adjacency matrix $\A_k$ of $F_k(G)$ that corresponds to the subscripts $X, X_{-x}^{+y}$ is $a_{xy}=1$  for every $x \in X$, $y \not\in X$ and $\{x,y\}\in E$. In all the other cases, $a_{xy}=0$. The Rayleigh quotient of $\vecv\in[n]_k$ is
 %Using these facts, we obtain:
\begin{equation}
\label{lambda(v)}
\lambda(\vecv)=\frac{\vecv^{\top}\L_k\vecv}{\vecv^{\top}\vecv} =\frac{\vecv^{\top}\D\vecv-\vecv^{\top}\A\vecv}{\vecv^{\top}\vecv}= \frac{P_k(\vecv)-Q_k(\vecv)-R_k(\vecv)}{S_k(\vecv)},
\end{equation}
%\blue{How about
%\begin{align*}
%\frac{\vecv^{\top}\L_k\vecv}{\vecv^{\top}\vecv} &= \frac{\sum_{\{X,Y\}\in %E(F_k)}[\vecv(X)-\vecv(Y)]^2}{\sum_{X\in V(F_k)}\vecv(X)^2}\\
% &=\frac{\sum_{X\in V(F_k)}\sum_{x\in X}\sum_{y\not\in %X}a_{xy}[\vecv(X)-\vecv(X_{-x}^{+y})]^2}{\sum_{X\in V(F_k)}\vecv(X)^2}
%\end{align*}
%}
where $P_k(\vecv)-Q_k(\vecv) =\vecv^{\top}\D\vecv$,
$R_k(\vecv)=\vecv^{\top}\A\vecv$, and $S_k(\vecv)=\vecv^{\top}\vecv$. Thus, according to the above comments, these expressions are
\begin{align}
P_k(\vecv) &=  \sum_{X \in [n]_k} \vecv(X)^2\left(\sum_{x \in X}d_{x} \right),\label{Pk}\\
Q_k(\vecv) &=  \sum_{X \in [n]_k}\vecv(X)^2\left(\sum_{ x,y  \in X}a_{xy} \right),\label{Qk}\\
R_k(\vecv) &= \sum_{X \in [n]_k}\sum_{x\in X}\sum_{y \not \in X}  \left(a_{xy}\vecv(X)\vecv(X_{-x}^{+y}) \right),\label{Rk}\\
S_k(\vecv) &= \sum_{X \in [n]_k} \vecv(X)^2.\label{Sk}
\end{align}

Through these expressions, the idea is to obtain a relationship between
the Rayleigh quotients of $\vecv$, indexed by the ${n\choose k}$ vertices  of $F_k(G)$, and vectors $\vecv_z^{\ast}$, with $z\in {n}$, indexed by the ${n\choose k-1}$ vertices of $F_{k-1}(G)$ for $k>1$. Using the same notation as in the proof of Theorem \ref{TH:Fk-alpha}, we first consider the vector $\vecv_z$, for $z\in[n]$, indexed by the vertices of $H_z$, and with entries $\vecv_z(Z)=\vecv(Z^{+z})$, for $Z\in [n\setminus\{z\}]_{k-1}$.
Then, the vector $\vecv_z^{\ast}$ is obtained from the vector $\vecv_z$, of dimension ${n-1\choose k-1}$, by adding ${n-1\choose k-2}$ null entries. Thus, 
the dimension of $\vecv_z^{\ast}$ is ${n-1\choose k-1}+{n-1\choose k-2}={n\choose k-1}$ and, for every $Z\in [n]_{k-1}$,
$$
\vecv_z^{\ast}(Z)=
\left\{
\begin{array}{ll}
\vecv_z(Z) & \mbox{if $z\not\in Z$,}\\
0 & \mbox{if $z\in Z$.}
\end{array}
\right.
$$ 
%  Given $Z$ a collection of $(k-1)$ distinct indices from $[n]$,  we assign the value  $\underline{\vecv}_z^Z$ to the entry from $\underline{\vecv}_z$ that corresponds to the configuration in which we have placed the $k-1$ distinct tokens to the nodes given by the indices from $Z$. In this way,  we obtain the next key relationships:
% Rearranging the terms, we get the following alternative expressions.
Then, when computing the Rayleigh quotient of $\vecv_z^{\ast}$ in $F_{k-1}(G)$, the above formulas of $P_k(\vecv)$, $Q_k(\vecv)$, $R_k(\vecv)$, and $S_k(\vecv)$ in \eqref{Pk}--\eqref{Sk} become
\begin{align}
P_{k-1}(\vecv_z^{\ast}) &=  \sum_{Z \in [n]_{k-1}} \vecv_z^{\ast}(Z)^2\left(\sum_{x \in Z}d_{x} \right),\label{Pk-1}\\
Q_{k-1}(\vecv_z^{\ast}) &=  \sum_{Z \in [n]_{k-1}}\vecv_z^{\ast}(Z)^2\left(\sum_{ x,y  \in Z}a_{xy} \right),\label{Qk-1}\\
R_{k-1}(\vecv_z^{\ast}) &= \sum_{Z \in [n]_{k-1}}\sum_{x\in Z}\sum_{y \not \in Z}  \left(a_{xy}\vecv_z^{\ast}(Z)\vecv_z^{\ast}(Z_{-x}^{+y}) \right),\label{Rk-1}\\
S_{k-1}(\vecv_z^{\ast}) &= \sum_{Z \in [n]_{k-1}} \vecv_z^{\ast}(Z)^2.\label{Sk-1}
\end{align}

Rearranging the terms of $P_k(\vecv)$, $Q_k(\vecv)$, $R_k(\vecv)$, and $S_k(\vecv)$, we obtain the relationships between them and $P_{k-1}(\vecv_z^{\ast})$, $Q_{k-1}(\vecv_z^{\ast})$, $R_{k-1}(\vecv_z^{\ast})$, and $S_{k-1}(\vecv_z^{\ast})$ in \eqref{Pk-1}--\eqref{Sk-1}, as follows:
\begin{align*}
P_k(\vecv) &= \frac{1}{k-1}\sum_{X \in [n]_{k}}  \vecv(X)^2\left(\sum_{z \in X} \sum_{x \in X_{-z} }d_{x} \right)=\frac{1}{k-1}\sum_{X \in [n]_{k}}  \sum_{z\in X} \vecv(X)^2\left( \sum_{x \in X_{-z}}d_{x} \right)\\
 &=\frac{1}{k-1}\sum_{z\in [n]} \sum_{Z \in [n\setminus \{z\}]_{k-1}} \vecv(Z^{+z})^2\left( \sum_{x \in Z }d_{x} \right)=
 \frac{1}{k-1}\sum_{z\in [n]}\sum_{Z \in [n]_{k-1}} \vecv_z^{\ast}(Z)^2\left(\sum_{x \in Z }d_{x} \right)\\
  &=\frac{1}{k-1} \sum_{z \in [n] } P_{k-1}(\vecv_z^{\ast}).
 \end{align*}
 For the first equality, notice that, for every given $X\in{[n]}_k$ and $x\in X$, the term $d_x$ appears once in $\sum_{x\in X}d_x$ but $k-1$ times in $\sum_{z\in X}\sum_{z\in X_{-z}}d_x$, and similarly for the other cases below. 
\begin{align*}
Q_k(\vecv) & = \frac{1}{k-2}\sum_{X \in [n]_{k}}  \vecv(X)^2\left(\sum_{z \in X} \sum_{x,y \in X_{-z} }a_{xy} \right)=\frac{1}{k-2}\sum_{X \in [n]_{k}}  \sum_{z\in X} \vecv(X)^2\left( \sum_{x,y \in X_{-z} }a_{xy} \right)\\
 &=\frac{1}{k-2}\sum_{z\in [n]} \sum_{Z \in [n\setminus \left\{z\right\}]_{k-1}} \vecv(Z^{+z} )^2\left( \sum_{x,y \in Z }a_{xy} \right)\\
 &=\frac{1}{k-2}\sum_{z\in [n]} \sum_{Z \in [n]_{k-1}} \vecv_z^{\ast}(Z)^2\left( \sum_{x,y \in Z }a_{xy} \right)
 = \frac{1}{k-2} \sum_{z \in [n] } Q_{k-1}(\vecv_z^{\ast})\quad \mbox{for $k>2$}.\\
 R_k(\vecv)&= \frac{1}{k-1} \sum_{X \in [n]_k}\sum_{z\in X}\sum_{x\in X \setminus \left\{z\right\} }\sum_{y \not \in X}  \left(a_{xy}\vecv(X)\vecv(X_{-x}^{+y} ) \right)\\
 &=\frac{1}{k-1} \sum_{z\in [n]} \sum_{Z \in [n \setminus \left\{z\right\}]_{k-1}}\sum_{x\in Z }\sum_{y \not \in Z^{+z} }  \left(a_{xy}\vecv(Z^{+z} )\vecv(Z_{-x}^{+y}) \right)
 \end{align*}
 \begin{align*}
 &=\frac{1}{k-1} \sum_{z\in [n]} \sum_{Z \in [n]_{k-1}}\sum_{x\in Z }\sum_{y \not \in Z }  \left(a_{xy}\vecv_x^{\ast}(Z)\vecv_z^{\ast} (Z_{-x}^{+y} ) \right)= \frac{1}{k-1} \sum_{z \in [n] } R_{k-1}(\vecv_z^{\ast}).\\
S_k(\vecv) &= \frac{1}{k} \sum_{z \in [n]} \sum_{Z \in [n \setminus \left\{z\right\} ]_{k-1}}\vecv(Z^{+z}  )^2 =\frac{1}{k} \sum_{z \in [n]} \sum_{Z \in [n \setminus \left\{z\right\} ]_{k-1}}\vecv_z(Z)^2=\frac{1}{k} \sum_{z \in [n]} \sum_{Z \in [n]_{k-1}}\vecv_z^{\ast}(Z)^2\\
	&= \frac{1}{k} \sum_{z\in [n]} S_{k-1}(\vecv_z^{\ast}).
\end{align*}

Putting these values into \eqref{lambda(v)}, we obtain
\begin{align}
\lambda(\vecv) &=
	%\frac{\vecv^{\top}L_k\vecv}{\vecv^{\top}\vecv} =
\frac{k}{k-1} \frac{\displaystyle\sum_{z\in [n]}\Big[(P_{k-1}(\vecv_z^{\ast})-Q_{k-1}(\vecv_z^{\ast})-R_{k-1}(\vecv_z^{\ast}))-\frac{1}{(k-2)}Q_{k-1}(\vecv_z^{\ast})\Big]}{\displaystyle\sum_{z \in [n]}S_{k-1}(\vecv_z^{\ast})}\nonumber\\
&= \frac{k}{k-1} \frac{\displaystyle\sum_{z\in [n]}\Big[(\vecv_z^{\ast})^{\top}\L_{k-1}\vecv_z^{\ast}-\frac{1}{(k-2)}Q_{k-1}(\vecv_z^{\ast})\Big]}{\displaystyle\sum_{z \in [n]}(\vecv_z^{\ast})^{\top}\vecv_z^{\ast}},\quad\text{ if $k>2$,} \label{k>2}\\
	\lambda(\vecv) &=	2\, \frac{\displaystyle\sum_{z\in [n]}\Big[(P_{1}(\vecv_z^{\ast})-R_{1}(\vecv_z^{\ast}))\Big]}{\displaystyle\sum_{z \in [n]}S_{1}(\vecv_z^{\ast})}-\frac{Q_{2}(\vecv)}{S_2(\vecv)},\quad \text{ if $k=2$.} \label{k=2}
\end{align}

Now suppose that $\lambda =\lambda(\vecv)= \alpha(F_k(G))$. From the hypothesis and \eqref{eq:inclusion}, $\alpha(F_k)$ is not an eigenvalue of $G$. Then, by Lemma \ref{lem:1}$(ii)$-$(iii)$, $\B^{\top}\vecv = \vec0$, and  both $\vecv_z$ and $\vecv_z^{\ast}$ are orthogonal to $\1$ (that is, they are embeddings of $F_{k-1}(G-z)$ and $F_{k-1}(G)$, respectively). Then,
$$
\lambda(\vecv_z^{\ast})=\frac{P_{k-1}(\vecv_z^{\ast})-Q_{k-1}(\vecv_z^{\ast})-R_{k-1}(\vecv_z^{\ast})}{S_{k-1}(\vecv_z^{\ast})} \geq \alpha(F_{k-1}(G)),
$$
or $(\vecv_z^{\ast})^{\top}\L_{k-1}\vecv_z^{\ast}\ge \alpha(F_{k-1}(G))(\vecv_z^{\ast})^{\top}\vecv_z^{\ast}$ for every $z\in[n]$, 
implying that
	\begin{equation} 
\label{bound-a(Fk-1)}
\frac{\displaystyle\sum_{z\in [n]}[(\vecv_z^{\ast})^{\top}\L_{k-1}\vecv_z^{\ast}]}
{\displaystyle\sum_{z \in [n]}[(\vecv_z^{\ast})^{\top}\vecv_z^{\ast}]} \geq \alpha(F_{k-1}(G)).
\end{equation}
Furthermore, since $a_{ij} \in \{0,1\}$ and $G$ has maximum degree $\Delta(G)$, we have 
\begin{equation}
\label{Qk-Sk}
Q_k(\vecw) \leq k\cdot\min\{k-1,\Delta(G)\}\cdot S_k(\vecw)
\end{equation}
Combining all the above results, from \eqref{k>2} to \eqref{Qk-Sk}, we obtain the claimed inequalities. Namely, if  $k \geq 3$,
\begin{align*}
\lambda =\alpha(F_k(G)) & \geq \frac{k}{k-1}\alpha(F_{k-1}(G))- \frac{k}{(k-1)(k-2)} \frac{\displaystyle\sum_{z\in [n]}Q_{k-1}(\vecv_z^{\ast})}{\displaystyle\sum_{z\in [n]}S_{k-1}(\vecv_z^{\ast})}\\
 &\geq \frac{k}{k-1} \alpha(F_{k-1}(G))-\frac{k}{k-2}\min\{k-2,\Delta(G)\},
\end{align*}
and, if $k = 2$,
\begin{align*}
	\lambda = \alpha(F_2(G)) &=2\, \frac{\displaystyle\sum_{z\in [n]}[(P_{1}(\vecv_z^{\ast})-R_{1}(\vecv_z^{\ast}))]}{\displaystyle\sum_{z \in [n]}S_{1}(\vecv_z^{\ast})}-\frac{Q_{2}(\vecv)}{S_2(\vecv)}   \geq 2 \alpha(G)- 2.
	\end{align*}
	This completes the proof.
\end{proof}

Let us now see that Theorem \ref{TH:Fk-alpha(b)}, with $\lambda=\alpha(F_k(G))$, is a consequence of  \eqref{eq:arnau}.  First,
notice that, in \eqref{eq:arnau}, we can replace $\min\{k-2,\Delta(G)\}$ by $k-2$ to have
\begin{equation}
\label{eq:arnau(b)}
\alpha_k(F_k(G))\ge \frac{k}{k-1}\alpha(F_{k-1}(G))-k.
\end{equation} 
Then, we use induction on $k\geq 2$ to show that 
\begin{equation}
\label{eq:Fk-alpha(bb)}
\alpha(F_k(G)) \geq k\alpha(G)-k(k-1). 
\end{equation} 
First, this is true for $k=2$, the base case. Now assume that, for some $k> 2$, \eqref{eq:Fk-alpha(bb)} holds. Then, using \eqref{eq:arnau(b)},
\begin{align*}
\alpha(F_{k+1}(G)) &\geq \frac{k+1}{k} \alpha(F_{k}(G)) - (k+1) \geq \frac{k+1}{k}( k\alpha(G)-k(k-1))-(k+1)\\
 &= (k+1) \alpha(G)-(k^2-1)-(k+1) = (k+1) \alpha(G)-(k+1)k.
 \end{align*}
 Alternatively, the same conclusion is reached by iterating \eqref{eq:arnau(b)}
 until the occurrence of $F_1(G)=G$.
 
Moreover, when $\Delta(G)$ is much smaller than $k-2$, the same method
 leads to the following result.
 \begin{corollary}
 If $k>2$, then
 \begin{equation}
 \label{eq:Fk-vs-G-Delta}
 \alpha(F_k(G))\ge k\, \alpha(G)-k\,\Delta(G)\left(1+\sum_{r=\Delta(G)}^{k-2}\frac{1}{r}\right)\ge  k\, \alpha(G)- k\Delta(G)\left(1+\ln\frac{k-1}{\Delta(G)}\right).
\end{equation}
 \end{corollary}

\end{document}